\documentclass[a4paper,10pt]{amsart}
\usepackage[utf8]{inputenc}
\usepackage{amsxtra}
\usepackage{amsopn}
\usepackage{amsmath,amsthm,amssymb}
\usepackage{amscd}
\usepackage{color}
\usepackage{amsfonts}
\usepackage{latexsym}
\usepackage{verbatim}

\usepackage{color}

\usepackage{verbatim}
\usepackage{comment}
\theoremstyle{plain}
\newtheorem{theorem}{Theorem}[section]
\newtheorem*{theorem*}{Theorem}

\newtheorem{lemma}[theorem]{Lemma}
\newtheorem{prop}[theorem]{Proposition}

\newtheorem{rem}[theorem]{Remark}

\newtheorem*{mt*}{Main Theorem}


\newcommand\C{{\mathbb C}}

\newcommand\Z{{\mathbb Z}}

\newcommand{\del}{\partial}
\newcommand{\delbar}{\overline{\del}}

\DeclareMathOperator{\Span}{Span}
\DeclareMathOperator{\id}{id}

\DeclareMathOperator{\im}{Im}
\DeclareMathSymbol{\Finv} {\mathord}{AMSb}{"60}

\title[Aeppli-Bott-Chern Massey products on non-K\"ahler solvmanifolds]{Aeppli-Bott-Chern Massey products on non-K\"ahler solvmanifolds}
\author[Nunzia Cesarino and Adriano Tomassini]{Nunzia Cesarino and Adriano Tomassini}
\address{Dipartimento di Scienze Matematiche, Fisiche e Informatiche\\
Unit\`{a} di Matematica e Informatica,
Universit\`{a} degli Studi di Parma\\
Parco Area delle Scienze 53/A, 43124 \\
Parma, Italy}
\email{nunzia.cesarino@unipr.it}
\email{adriano.tomassini@unipr.it}
\thanks{\newline 
The second author is partially supported by the Project PRIN 2022 ``Real and Complex Manifolds: Geometry and Holomorphic Dynamics 2022AP8HZ9'' and by GNSAGA of INdAM}

\begin{document}

\begin{abstract}
In this paper, we present explicit computations of non-trivial triple $ABC$–Massey products on non-K\"ahler solvmanifolds endowed with an invariant complex structure. We prove that the {\em Bigalke-Rollenske manifold}, the {\em generalized Nakamura manifolds} satisfying some suitable assumptions and 
compact quotients of the solvable Lie group $\mathbb{C}^{2n}\ltimes_{\rho} \mathbb{C}^{2m}$ have non-vanishing triple $ABC$-Massey products. Furthermore, such manifolds have no astheno-K\"ahler metric.
\end{abstract}
\maketitle

\section{Introduction}
A useful tool to show that a given smooth manifold has no Kähler structure is to construct a non-vanishing triple Massey product, a cohomology class living in a suitable quotient of the de Rham cohomology.\\
\indent
In the framework of rational homotopy theory, Sullivan \cite{Sul77} defines a commutative differential graded algebra (cdga) as \emph{formal} if it is connected to its (de Rham) cohomology by a chain of multiplicative quasi-isomorphisms.\\
\indent As a consequence, for simply connected formal spaces, the rational homotopy groups $\pi_{\ge 2} \otimes \mathbb{Q}$ can be recovered directly from the cohomology ring, and Massey products, which naturally arise as higher-order cohomological operations in the quasi-isomorphism class of a cdga,  constitute the main obstructions to formality.
%vanishing in the formal case. 
Furthermore, he observed \cite{Sul75} that
“there are topological obstructions for $M$ to admit a metric in which the product of
harmonic forms is still harmonic”: for instance, the rational formality of the closed
oriented manifold and the vanishing of every Massey product. 
%If a manifold $M$ admits such a metric, then it is formal;
Consequently, non-formal manifolds cannot admit a metric for which the wedge product of harmonic forms is still harmonic. Motivated by this fact, Kotschick introduced the notion of \emph{geometrically formal manifold} \cite[Definition 1]{Kot}, i.e., a closed manifold admitting a Riemannian metric satisfying such a property.

On a complex manifold, the decomposition $d=\del+\delbar$ gives rise to suitable complex cohomologies, namely the {\em Bott-Chern} and {\em Aeppli} cohomology groups. It turns out that when $M$ is closed and carries a K\"ahler metric, then all such cohomologies are isomorphic to the {\em Dolbeault} one. Therefore,  Bott-Chern and Aeppli cohomologies provide additional data that can be useful for the study of compact non-K\"ahler manifolds. 

Inspired by Kotschick, {\em geometrically Bott--Chern and Aeppli, ABC-geometrically} {\em  formal} metrics \cite{AT-1, ASfe, MS} are introduced. These metrics are indeed defined through the Bott-Chern and Aeppli Laplacians (see Section \ref{preliminari} for the definitions) and arise naturally in pluripotential homotopy theory \cite{MS}. 
\vskip5truecm
In \cite[Definition 2.1]{AT-1}, the authors, motivated by the search for complex invariants and obstructions to the existence of geometrically Bott-Chern formal metrics on a compact complex manifold, defined and studied the \emph{triple Aeppli--Bott--Chern--Massey products} (briefly, $ABC$-{\em Massey products}). Recently, Milivojevi\'c and Stelzig in \cite{MS} gave the definition of higher Aeppli–Bott–Chern–Massey products, which extend the  triple $ABC$-Massey products. For other results on triple $ABC$-Massey products, we refer to \cite{TarTom, STc, PSZ24} and the references therein.
\vskip.1truecm
The aim of this paper is to present explicit computations of non-trivial triple $ABC$–Massey products on solvmanifolds, namely compact quotients of a simply connected solvable Lie groups by a lattice, endowed with an invariant complex structure. The cohomological and metric properties of such manifolds have been extensively studied by many authors, since solvmanifolds provide a source of concrete 
examples admitting further structures--e.g., complex structures, symplectic structures, special metrics as balanced, strong K\"ahler with torsion, and astheno-K\"ahler metrics--on which explicit cohomological computations can be carried out at Lie algebra level (see \cite{AngKas, AngKas2, Kas, Has05, Has06, R, ASfe, FV, AN, STa} and the references therein). \newline 
More precisely, after recalling some preliminaries of complex and Hermitian geometry and fixing some notation (see Section \ref{preliminari}), in Section \ref{main} we prove that the {\em Bigalke-Rollenske manifold} 
$(M^{4n-2},J)$ \cite{BigRol} does not admit any geometrically-Bott-Chern-formal metric and has no astheno-K\"ahler metric (see Theorem \ref{BRoll}). \newline 
In Section \ref{NakamuraM} we focus on the {\em generalized Nakamura manifolds} $N_{M,P,\tau}$, introduced in \cite{CattTom} as an extension of the {\em Nakamura manifolds}  \cite{Nak}: we observe that if $N_{M,P,\tau}$ satisfies the $\del\delbar$-Lemma, then it admits a geometrically-Bott-Chern-formal metric (see Theorem \ref{del-delbar}). Then, under suitable assumptions, we prove that $N_{M,P,\tau}$ has no geometrically Bott–Chern formal metrics (see Theorem \ref{gNakamura}) and has no astheno-K\"ahler metrics (see Proposition \ref{no-astheno}). It has to be remarked that, if a compact complex manifold satisfies the $\partial\bar{\partial}$-Lemma, then it is formal and consequently all Massey products vanish. In contrast, there exist compact complex manifolds satisfying the $\partial\bar{\partial}$-Lemma admitting non-trivial $ABC$-Massey products (see \cite{SfeTom1}).
However, in \cite[Corollary 5.11]{STb}, it is shown that the non-vanishing of $ABC$--Massey products detects the non-K\"ahlerianity of solvmanifolds. \newline Section \ref{LT} is devoted to the construction of non-trivial triple $ABC$-Massey products on compact quotients $N_{\underline{\mu},P}$ of the Lie group $\mathbb{C}^{2n}\ltimes_{\rho} \mathbb{C}^{2m}$. Such manifolds have been recently introduced in \cite{LuTom} and provide families of compact complex manifolds with no K\"ahler structure and carrying symplectic structures satisfying the Hard Lefschetz Condition. Then we prove that $N_{\underline{\mu},P}$ has no geometrically-Bott-Chern-formal metric neither astheno-K\"ahler metric (see Theorem \ref{Lu1} and Proposition \ref{Lu2}).
 The proofs of our results involve tools from Hodge theory for the Bott-Chern and Aeppli Laplacians \cite{Schw} as well as the computation of suitable spaces of harmonic forms, which allow us to construct triple $ABC$-Massey products and to show that such products do not vanish.
\vskip.3truecm
{\em \underline{Acknowledgments:}} The authors would like to thank Lorenzo Sillari and Ettore Lo Giudice for useful comments and remarks. Many thanks are also due to Jonas Stelzig for bringing to their attention the reference \cite{trivial-bundle}.

\section{Preliminaries}\label{preliminari}
\label{preliminaries}
In this Section we recall some facts about complex and Hermitian manifolds and fix some notation.
Let $M$ be a smooth compact manifold of dimension $2n$ and let $J$ be a complex structure on $M$, namely a $(1,1)$-tensor on $M$ such that $J^2=-\id$ satisfying the integrability condition
$$
N_J(U,V)=[JU,JV]-[U,V]-J[JU,V]-J[U,JV]=0,
$$ 
for every pair of vector fields $U,V$ on $M$.
Then, $J$ induces on the space of complex forms $\mathcal{A}^{\bullet}_{\C}(M)$ a natural bigrading, namely $\mathcal{A}^{\bullet}_{\C}(M)=\bigoplus_{p+q=\bullet}{\mathcal A}^{p,q}(M)$, where ${\mathcal A}^{p,q}(M)$ denotes the space of $(p,q)$-forms on $M$. Accordingly, the exterior derivative $d$ acts on $(p,q)$-forms as
$$
d:{\mathcal A}^{p,q}(M)\to  {\mathcal A}^{p+1,q}(M)\oplus {\mathcal A}^{p,q+1}(M)
$$
so that $d$ decomposes into
$$
d=\del+\delbar,
$$
and the differential operators $\del,\delbar$ satisfy the following
$$
\del^2=0,\quad \delbar^2=0,\qquad \del\delbar+\delbar\del=0.
$$
The {\em Dolbeault}, {\em Bott-Chern} and {\em Aeppli} cohomology groups of the complex manifold $(M,J)$ are defined, respectively, as
$$
H_{\delbar}^{\bullet,\bullet}(M):=\frac{\ker\delbar}{\im\delbar}\cap \mathcal{A}^{\bullet,\bullet}(M),
$$
$$
H_{BC}^{\bullet,\bullet}(M):=\frac{\ker \del \cap \ker\delbar}{\im \del\delbar}\cap \mathcal{A}^{\bullet,\bullet}(M),
$$
$$
H_{A}^{\bullet,\bullet}(M):=\frac{\ker \del\delbar}{\im \del+\im \delbar}\cap \mathcal{A}^{\bullet,\bullet}(M).
$$
Let $g$ be a Hermitian metric on a complex manifold $(M,J)$. According to \cite{Schw}, the \emph{Bott-Chern Laplacian} and the \emph{Aeppli Laplacian} are respectively defined as
\begin{eqnarray*}
\Delta_{BC}&:=& \del\delbar\delbar^*\del^* +\delbar^*\del^*\del\delbar + \delbar^*\del \del^*\delbar + \del^*\delbar\delbar^*\del +\delbar^*\delbar +\del^*\del,\\
\Delta_{A}&:=& \del\del^*+\delbar\delbar^*+ \delbar^*\del^*\del\delbar +\del\delbar\delbar^*\del^*+\del\delbar^*\delbar\del^*+\delbar\del^*\del\delbar^*.
\end{eqnarray*}
Both $\Delta_{BC}$ and $\Delta_{A}$ are fourth order self-adjoint elliptic differential operators and, consequently, the following isomorphisms of vector spaces hold
$$
\mathcal{H}_{\Delta_{BC}}^{\bullet,\bullet}(M):=\ker\Delta_{BC}\vert_{\mathcal{A}^{\bullet,\bullet}(M)}\cong H^{\bullet,\bullet}_{BC}(M),$$  
$$ \mathcal{H}_{\Delta_{A}}^{\bullet,\bullet}(M):=\ker\Delta_{A}\vert_{\mathcal{A}^{\bullet,\bullet}(M)}\cong H^{\bullet,\bullet}_{A}(M).
$$
It turns out that a $(p,q)$-form $\alpha$ is \emph{Bott-Chern harmonic} (respectively, \emph{Aeppli harmonic}) with respect to a fixed Hermitian metric $g$ on $(M,J)$ if and only if the following conditions hold 
$$
\del\alpha=0,\qquad \delbar\alpha=0 \qquad \del\delbar\ast\alpha=0,
$$
respectively,
$$
\del\delbar\alpha=0, \qquad \del\ast\alpha=0, \qquad \delbar\ast\alpha=0,
$$
where $\ast: \mathcal{A}^{p,q}(M)\to \mathcal{A}^{n-p,n-q}(M)$ is the $\mathbb{C}$-antilinear Hodge $\ast$ operator with respect to the Hermitian metric $g$.\\
Following Kotschick \cite{Kot}, in \cite[Definition 1.1]{AT-1} the authors define a Hermitian metric $g$ on a compact complex manifold $(M,J)$ to be \emph{geometrically-Bott-Chern-formal} (shortly \emph{geometrically-$BC$-formal}) if the space of Bott-Chern harmonic forms $\mathcal{H}_{\Delta_{BC}}^{\bullet,\bullet}(M)$
is endowed with an algebra structure induced by the wedge product.\\
\indent We now recall the construction of \emph{triple $ABC$-Massey products} (\cite[Definition 2.1]{AT-1}).
Let $\mathfrak {a}_{12}=[\alpha_{12}]\in H_{BC}^{p,q}(M)$, $\mathfrak {a}_{23}=[\alpha_{23}]\in H_{BC}^{r,s}(M)$, and $\mathfrak {a}_{34}=[\alpha_{34}]\in H_{BC}^{u,v}(M)$ such that 
$$
\mathfrak {a}_{12}\cup\mathfrak {a}_{23}=0,\qquad \mathfrak {a}_{23}\cup\mathfrak {a}_{34}=0.
$$
Let 
$$
(-1)^{p+q}\alpha_{12}\wedge\alpha_{23}=\del\delbar f_{13}, \qquad (-1)^{r+s}\alpha_{23}\wedge\alpha_{34}=\del\delbar f_{24}.
$$
Then, the \emph{triple $ABC$-Massey product} $\langle\mathfrak {a}_{12},\mathfrak {a}_{23},\mathfrak {a}_{34}\rangle_{ABC}$ is defined as
\begin{eqnarray*}
\langle\mathfrak {a}_{12},\mathfrak {a}_{23},\mathfrak {a}_{34}\rangle_{ABC}&=&[(-1)^{p+q}\alpha_{12}\wedge f_{24}-(-1)^{r+s}f_{13}\wedge\alpha_{34}]\\[5pt]
&\in& \frac{H_A^{p+r+u-1,q+s+v-1}(M)}{\mathfrak {a}_{12}\cup H_A^{r+u-1,s+v-1}(M)+ H_A^{p+r-1,q+s-1}(M)\cup\mathfrak {a}_{34}}.
\end{eqnarray*}
Triple $ABC$-Massey products provide an obstruction to the existence of geometrically-Bott-Chern-formal metrics on a compact complex manifold \cite[Theorem 2.4]{AT-1}, that is, if a compact complex manifold admits a geometrically-Bott–Chern-formal metric, then the triple Aeppli–Bott–Chern–Massey products vanish.\\
\indent We next recall the definition of \emph{astheno-Kähler} metrics. Given a Hermitian metric $g$ on a complex manifold $(M,J)$, the associated fundamental form $\omega$ is defined by
$$\omega(U,V)=g(JU,V),\qquad \forall\, U,V \in \Gamma(TM).$$
A Hermitian metric $g$ on $(M,J)$ is said to be \emph{astheno-K\"{a}hler} if its fundamental form $\omega$ satisfies the condition $$\del\delbar\omega^{n-2}=0.\\$$
We denote by $\Omega=\frac{\omega^{n}}{n!}$ the standard volume form associated to the Hermitian metric $g$. 
Let $x\in M$ and let $\psi\in\Lambda^{p,0}_{\mathbb{R}}(T_{x}^{\ast}M)$. Then, $\psi$ is said to be {\em simple} if
    $$\psi=\psi^{1}\wedge\dots\wedge\psi^{p},$$ where $\psi^{i}\in\Lambda^{1,0}(T^{\ast}_{x}M\otimes\mathbb{C})$ for $i=1,\dots, p$. Let $\psi\in\Lambda^{n,n}_{\mathbb{R}}(T_{x}^{\ast}M)$. Then, $\psi$ is said to be {\em strictly positive} if 
$$
\psi=a\Omega_x,
$$
for $a>0$. Let  $\psi\in\Lambda^{p,p}_{\mathbb{R}}(T_{x}^{\ast}M)$. Then $\psi$ is said to be \emph{transverse} if for any given simple $\beta\in \Lambda^{n-p,0}(T^{\ast}_{x}M\otimes\mathbb{C})$ with $\beta\ne 0$, 
   $$\sigma_{n-p}\,\psi\wedge\beta\wedge\overline{\beta}$$
is strictly positive. It turns out that (see \cite[Corollary 1.10]{HK}), if $\omega$ is the fundamental form of a Hermitian metric $g$ on $M$, then $\omega^p$ is transverse for $1\leq p \leq n$.

\indent We now turn to some basic facts of complex geometry on Lie groups. Let $G$ be a real $2n$-dimensional, connected and simply connected Lie group with Lie algebra $\mathfrak{g}$. If $\Gamma$ is a lattice in $G$, we may consider the compact manifold given by the quotient $\Gamma\backslash G$. In particular, if $G$ is nilpotent we call $\Gamma\backslash G$ a \emph{nilmanifold}; if $G$ is solvable we call $\Gamma\backslash G$ a \emph{solvmanifold}.\\
An almost complex structure $J$ on a real $2n$-dimensional Lie algebra $\mathfrak{g}$ can be defined by assigning an $n$-dimensional complex subspace   $$\mathfrak{g}_{1,0}
\subset\mathfrak{g}_{\mathbb{C}}:
=\mathfrak{g}\otimes\mathbb{C}$$
of the complexified Lie algebra $\mathfrak{g}_{\mathbb{C}}$ such that 
$$\mathfrak{g}_{1,0}\cap\,\overline{\mathfrak{g}_{1,0}}
=\{0\},$$
by setting $$ 
\begin{array}{ll}
J(v):=\sqrt{-1}\,  v &\qquad \text{for} \quad v\in\mathfrak{g}_{1,0},\\[5pt]
J(w):=-\sqrt{-1}\,  w & \qquad \text{for} \quad w\in\overline{\mathfrak{g}_{1,0}}. 
\end{array}
$$
The almost complex structure $J$ is integrable if 
$$ [\mathfrak{g}_{1,0},\mathfrak{g}_{1,0}]\subset\mathfrak{g}_{1,0}.
$$
Any left-invariant complex structure $J$ on  $\mathfrak{g}$ determines a left-invariant complex structure on $G$, and hence induces one on $\Gamma\backslash G$, which we also denote by $J$. Equivalently, if $$\mathfrak{g}^{1,0}=\Span_\C\langle \phi^1,\ldots,\phi^{n}\rangle\subset \mathfrak{g}_{\mathbb{C}}^{\ast},$$
then the integrability condition of the almost complex structure induced becomes equivalent to the following 
$$
d\phi^j\in\Lambda^{2,0}\mathfrak{g}_{\mathbb{C}}^{\ast}\oplus \Lambda^{1,1}\mathfrak{g}_{\mathbb{C}}^{\ast},\qquad j=1,\ldots, n.
$$
A left-invariant complex structure $J$ is called \emph{nilpotent} \cite{JNilp} if there exists a co-frame $\{\varphi^{1},\dots, \varphi^{n}\}$ of left-invariant $(1,0)$-forms such that
\begin{equation*}
    d\varphi^{i}=\Lambda^{2}\langle \varphi^{1}, \dots, \varphi^{i-1},\overline{\varphi}^{1}, \dots , \overline{\varphi}^{i-1}\rangle, \quad\quad i=1,\cdots n.
\end{equation*}
\indent Finally, we recall the definition of {\em Kodaira dimension}. Let $(M,J)$ be an $n$-dimensional complex manifold and let $\Lambda^{n,0}(M)$ denote the canonical bundle of $(M,J)$. Let $$P_r(M,J)=\dim\,H^0(M,(\Lambda^{n,0}(M))^{\otimes r})$$ the $r^{th}$ {plurigenus} of $(M,J)$. Then the {\em Kodaira dimension} of the complex manifold $(M,J)$ is defined as
$$
\kappa^J(M)=
\left\{
\begin{array}{ll}
-\infty & \hbox{\rm if}\,\, P_r(M,J)=0\,\, \hbox{\rm for all} \,\,r\geq 1,\\[10pt]
\limsup_{r\to\infty}\frac{\log P_r(M,J)}{\log r} & \hbox{\rm otherwise.}
\end{array}
\right.
$$
\section{Non-formality of Bigalke-Rollenske manifolds}\label{main}
We begin by recalling the construction of the {\em Bigalke-Rollenske manifold}, see  \cite{BigRol} for the definition and \cite{SfeTar} for the existence of special metrics on it. \\For $n\geq 2$, let $G_n$  be the real nilpotent subgroup of $GL(2n+2;\C)$, whose Lie algebra will be denoted by $\mathfrak{g}$, consisting of upper triangular matrices of the form
$$
A = 
\begin{pmatrix}
1 &0 & &\cdots &\cdots &\cdots &0& \bar{y}_1& w_1 \\
&1& 0 \cdots 0 &\bar{z}_1& -x_1&0& \cdots &0& w_2 \\
& & \ddots  & \ddots &\ddots&\ddots&\ddots& \vdots&\vdots  \\
& & &1& 0 \cdots 0&\bar{z}_{n-1} &-x_{n-1} &0& w_n \\
& & & & 1&0& \cdots &0& y_1 \\
& & & & &\ddots& & \vdots&\vdots\\
& & & & &          & 1 & 0 & y_n\\
& & & & &          &  & 1& z_1\\
& & & & &          &  & & 1\\
 \end{pmatrix},
$$
where $x_1,\ldots,x_{n-1}, y_1,\ldots,y_n,z_1,\ldots,z_{n-1},w_1,\ldots,w_n\in\C$. \newline 
Denoting by $$\Gamma=G_n\cap GL(2n+2;\Z[\sqrt{-1}])$$ the subgroup of $G_n$ whose matrices have entries in $\Z[\sqrt{-1}]$, it turns out that $\Gamma$ is a lattice in $G_n$ and, consequently, the quotient $$M^{4n-2}:=\Gamma\backslash G_n$$ is a compact nilmanifold of complex dimension $4n-2$. By setting
$$
\left\{
\begin{array}{ll}
\phi^1&=dx_1,\\
&\vdots\\
\phi^{n-1}&=dx_{n-1},\\
\phi^{n}&=dy_{1},\\
&\vdots\\
\phi^{2n-1}&=dy_{n},\\
\phi^{2n}&=dz_{1},\\
&\vdots\\
\phi^{3n-2}&=dz_{n-1},\\
\phi^{3n-1}&=dw_1-\overline{y}_1dz_{1},\\
\phi^{3n}&=dw_2-\overline{z}_1dy_{1}+x_{1}dy_2,\\
&\vdots\\
\phi^{4n-2}&=dw_n-\overline{z}_{n-1}dy_{n-1}+x_{n-1}dy_n,\\
\end{array}
\right.
$$
a direct computation shows that $\phi^{1},\ldots,\phi^{4n-2}$ are complex left-invariant $1$-forms on $G_n$ and 
that the following structure equations hold:
\begin{equation}
\left\{
\begin{array}{lll}\label{structure-equations-Big-Rol}
d\phi^j&=0, & 1\leq j\leq 3n-2,\\[10pt]
d\phi^{3n-1}&=\phi^{2n}\wedge\overline{\phi}^{n}, &\\[10pt]
d\phi^{j}&=\phi^{j-3n+1}\wedge\phi^{j-2n+1}+\phi^{j-2n}\wedge\overline{\phi}^{j-n},& 3n\leq j\leq 4n-2.\\ 
\end{array}
\right.
\end{equation}
Therefore, by setting 
$$
\mathfrak{g}^{1,0}:=\Span_\C\langle \phi^1,\ldots,\phi^{4n-2}\rangle,
$$
we obtain a left-invariant almost complex structure $J$ on $M^{4n-2}$ induced by the decomposition 
$$\mathfrak{g}^*_\C=\mathfrak{g}^{1,0}\oplus
\overline{\mathfrak{g}^{1,0}},$$ 
and equations \eqref{structure-equations-Big-Rol} imply that $J$ is integrable. The complex nilmanifold $(M^{4n-2},J)$ is called the {\em Bigalke-Rollenske manifold}. It follows by the structure equations that the complex structure $J$ is nilpotent. We are ready to prove the following:
\begin{theorem}\label{BRoll}
Let $M=(M^{4n-2},J)$ be the Bigalke-Rollenske manifold. Then
\begin{enumerate}
\item[1)] $M$ does not admit any geometrically-Bott-Chern-formal metric;\\[5pt]
\item[2)] $M$ has no astheno-K\"ahler metrics;\\[5pt]
\item[3)] $\kappa^J(M)=0$.
\end{enumerate}
\end{theorem}
\begin{proof}
1) In order to show that $M$ does not carry any geometrically-Bott-Chern-formal metric, according to \cite[Theorem 2.4]{AT-1}, it is enough to construct a non-vanishing triple $ABC$-Massey product on it. To this purpose, let us take the following classes in Bott-Chern, defined respectively as
$$
\mathfrak{a}_{12}=[\alpha_{12}]=[\phi^n\wedge\overline{\phi}^n], \qquad \mathfrak{a}_{23}=[\alpha_{23}]=[\phi^{2n}\wedge\overline{\phi}^{2n}],\qquad
\mathfrak{a}_{34}=[\alpha_{34}]=[\phi^{2n}\wedge\overline{\phi}^{2n}].
$$
Then, 
$$\mathfrak{a}_{12}\cup \mathfrak{a}_{23}=0,\quad \mathfrak{a}_{23}\cup\mathfrak{a}_{34}=0,$$
since
$$
\alpha_{12}\wedge \alpha_{23}=\phi^n\wedge\overline{\phi}^n\wedge\phi^{2n}\wedge\overline{\phi}^{2n}=\del\delbar(\phi^{3n-1}\wedge\overline{\phi}^{3n-1})\,,
\quad \alpha_{23}\wedge \alpha_{34}=0.
$$
Therefore, the following $ABC$-Massey product is well-defined:
\begin{eqnarray*}
\langle\mathfrak {a}_{12},\mathfrak {a}_{23},\mathfrak {a}_{34}\rangle_{ABC}&=&[-\phi^{3n-1}\wedge\overline{\phi}^{3n-1}\wedge\phi^{2n}\wedge\overline{\phi}^{2n}]\\[5pt]
&\in& \frac{H_A^{2,2}(M)}{[\phi^n\wedge\overline{\phi}^n]\cup H_A^{1,1}(M)+ H_A^{1,1}(M)\cup[\phi^{2n}\wedge\overline{\phi}^{2n}]}.
\end{eqnarray*}
We are going to show that
$$
[-\phi^{3n-1}\wedge\overline{\phi}^{3n-1}\wedge\phi^{2n}\wedge\overline{\phi}^{2n}]\notin
[\phi^n\wedge\overline{\phi}^n]\cup H_A^{1,1}(M)+ H_A^{1,1}(M)\cup[\phi^{2n}\wedge\overline{\phi}^{2n}],
$$
so that $\langle\mathfrak {a}_{12},\mathfrak {a}_{23},\mathfrak {a}_{34}\rangle_{ABC}$ is a non-vanishing triple $ABC$-Massey product on $M$.\\
By contradiction, suppose that there exist $\alpha^{1,1},\, \beta^{1,1}\in\mathcal{A}^{1,1}(M)$, $P^{1,2}_{1}\in\mathcal{A}^{1,2}(M)$ and $P_{2}^{2,1}\in\mathcal{A}^{2,1}(M)$  such that 
$$\
\del\delbar\,\alpha^{1,1}= \del\delbar\,\beta^{1,1}=0,
$$ 
and 
\begin{equation}\label{Aeppli}
   \phi^{2n}\wedge\phi^{3n-1}\wedge\overline{\phi}^{2n}\wedge\overline{\phi}^{3n-1}=\phi^{n}\wedge\overline{\phi}^{n}\wedge \alpha^{1,1}+\beta^{1,1}\wedge\phi^{2n}\wedge\overline{\phi}^{2n}+\del P_{1}^{1,2}+\delbar P_{2}^{2,1}.
\end{equation}
Let $g$
%$$
%g=\sum_{j=1}^{4n-2}\phi^{j}\otimes\overline{\phi}^{j},
%$$
be the invariant Hermitian metric on $M$ whose associated fundamental form is given by 
$$\omega=\frac{\sqrt{-1}}{2}\sum_{j=1}^{4n-2} \phi^{j}\wedge\overline{\phi}^{j},$$ and let $\Omega=\frac{\omega^{4n-2}}{(4n-2)!}$ denote the standard volume form.\\  Multiplying equation \eqref{Aeppli} by $\ast(\phi^{2n}\wedge\phi^{3n-1}\wedge\overline{\phi}^{2n}\wedge\overline{\phi}^{3n-1})$, we obtain 
\begin{equation*}
    C\,\Omega= \ast(\phi^{2n}\wedge\phi^{3n-1}\wedge\overline{\phi}^{2n}\wedge\overline{\phi}^{3n-1})\wedge (\beta^{1,1}\wedge\phi^{2n}\wedge\overline{\phi}^{2n}+\del P_{1}^{1,2}+\delbar P_{2}^{2,1})
\end{equation*}
for $$C=g(\phi^{2n}\wedge\phi^{3n-1}\wedge\overline{\phi}^{2n}\wedge\overline{\phi}^{3n-1}, \phi^{2n}\wedge\phi^{3n-1}\wedge\overline{\phi}^{2n}\wedge\overline{\phi}^{3n-1})\in\mathbb{C}\setminus\{0\}.$$
A direct computation shows that $\phi^{2n}\wedge\phi^{3n-1}\wedge\overline{\phi}^{2n}\wedge\overline{\phi}^{3n-1}$ is Aeppli-harmonic, so that 
$$d^{\ast}(\phi^{2n}\wedge\phi^{3n-1}\wedge\overline{\phi}^{2n}\wedge\overline{\phi}^{3n-1})=0,$$ and we can then rewrite the last equation as
\begin{equation}\label{Aeppli1}
    C\,\Omega=\ast(\phi^{2n}\wedge\phi^{3n-1}\wedge\overline{\phi}^{2n}\wedge\overline{\phi}^{3n-1})\wedge\beta^{1,1}\wedge\phi^{2n}\wedge\overline{\phi}^{2n}+d\, P
\end{equation}
where $$P=\ast(\phi^{2n}\wedge\phi^{3n-1}\wedge\overline{\phi}^{2n}\wedge\overline{\phi}^{3n-1})\wedge( P_{1}^{1,2}+ P_{2}^{2,1}).
$$
Since the left-invariant complex structure $J$ on $M$ is nilpotent, by applying \cite[Corollary 2.7]{AngKas} it turns out that the Aeppli and Bott-Chern cohomologies can be computed in terms of complex left-invariant forms on $M$. Therefore, 
$$
\ker\del\delbar\vert_{\Lambda^{1,1}\mathfrak{g}_{\mathbb{C}}^{\ast}}\supset\mathcal{H}_{\Delta_{A}}^{1,1}(\mathfrak{g}^*)=\ker\Delta_{A}\vert_{\Lambda^{1,1}\mathfrak{g}_{\mathbb{C}}^{\ast}}\cong H^{1,1}_{A}(M)
$$
Setting
\begin{equation*}
   \beta^{1,1}=\sum_{1\le i,\,j\le 4n-2}c_{i\overline{j}}\,\phi^{i}\wedge\overline{\phi}^{j},\quad\quad\text{with}\quad c_{i\overline{j}}\in\C, 
\end{equation*}
a straightforward calculation shows that if $\del\delbar \beta^{1,1}=0$, then the following condition holds
$$
c_{3n-1\,\overline{3n-1}}=0.
$$
Therefore, equation \eqref{Aeppli1} becomes
\begin{equation*}
    C\Omega= d\,P,
\end{equation*}
which leads to a contradiction, $C$ being a non-zero constant and $\Omega$ a volume form. \\
Thus, we showed that
\begin{equation*}
\langle\mathfrak {a}_{12},\mathfrak {a}_{23},\mathfrak {a}_{34}\rangle_{ABC}=[-\phi^{3n-1}\wedge\overline{\phi}^{3n-1}\wedge\phi^{2n}\wedge\overline{\phi}^{2n}]\ne 0,
\end{equation*}
which means that $M$ admits a non-vanishing $ABC$-Massey product.\\\\
2) We will produce a suitable $\del\delbar$-exact $(2,2)$-form. Then, applying a standard argument, which we will repeat for the sake of completeness (see e.g., \cite[Lemma 3.5]{STa}), we will conclude that there are no astheno-K\"ahler metrics on $M$. A direct computation using the structure equations \eqref{structure-equations-Big-Rol} yields
$$\del\delbar(-\phi^{3n-1}\wedge\overline{\phi}^{3n-1})=\phi^{n}\wedge\phi^{2n}\wedge\overline{\phi}^{n}\wedge\overline{\phi}^{2n}.$$
If $\omega$ were the fundamental form of an astheno-K\"ahler metric on $M$, we should get
\begin{eqnarray*}
0&<& \int_{M}\omega^{4n-4}\wedge\phi^{n}\wedge\phi^{2n}\wedge\overline{\phi}^{n}\wedge\overline{\phi}^{2n}\\
&=& -\int_{M}\omega^{4n-4}\wedge \del\delbar(\phi^{3n-1}\wedge\overline{\phi}^{3n-1})\\
&=&-\int_{M}\del\delbar(\omega^{4n-4})\wedge \phi^{3n-1}\wedge\overline{\phi}^{3n-1}=0.
\end{eqnarray*}
This is absurd.
%On the other hand, applying Stokes' theorem yields
%\begin{equation*}
 %  \int_{M}\omega^{4n-4}\wedge\del\delbar(\phi^{3n-1}\wedge\overline{\phi}^{3n-1})= -\int_{M}d\omega^{4n-4}\wedge\delbar(\phi^{3n-1}\wedge\overline{\phi}^{3n-1})=0
%\end{equation*}
%due to bidegree reasons. This contradicts the previous inequality.\\\\
\vskip.3truecm\noindent
3) It is a consequence of the triviality of the canonical bundle on nilmanifolds (see, e.g., \cite{trivial-bundle}). However, we will provide an explicit proof. \\
Let $\eta\in\mathcal{A}^{n,0}(M)$, then we can write
$$\eta=f \phi^{1}\wedge\dots\wedge\phi^{4n-2},\quad\quad\text{with}\quad f\in C^{\infty}(M,\mathbb{C}).$$
A direct calculation shows that
\begin{equation*}
    \delbar(\phi^{1}\wedge\dots\wedge\phi^{4n-2})=0.
\end{equation*}
Therefore,
\begin{equation*}
    0=\delbar\eta=(\delbar\,f)\wedge \phi^{1}\wedge\dots\wedge\phi^{4n-2}.
\end{equation*}
Since the wedge product \(\phi^{1} \wedge \dots \wedge \phi^{4n-2}\) is non-vanishing, it follows that
\[
\bar{\partial} f = 0,
\]
i.e., \(f\) is a holomorphic function on the compact complex manifold \(M\). As a consequence, $f$ must be constant and, therefore, $P_1(M,J)=1$. Similar computations give $P_r(M,J)=1$. 
Indeed, by induction it is immediate to see that $\delbar((\phi^1\wedge\cdots\wedge\phi^{4n-2})^{\otimes r})=0$ for every $r\geq 1$, so that the condition 
$$
\delbar(f(\phi^1\wedge\cdots\wedge\phi^{4n-2})^{\otimes r})=0
$$ 
is again equivalent to $\delbar f=0$. \\ Thus, $P_r(M,J) =1$ for every $r\geq 1$, and hence $\kappa^J(M) = 0$.
\end{proof}

%\section{Explicit constructions}\label{section:examples}

\section{Non-formality of generalized Nakamura manifolds}\label{NakamuraM}
We recall the construction of the {\em generalized Nakamura manifolds}, provided in \cite{CattTom}, which consists of a family of compact solvmanifolds of completely solvable type that are not holomorphically parallelizable.\\
Let $M\in SL(n;\mathbb{Z})$ be a diagonalizable matrix over $\mathbb{R}$ with positive eigenvalues $\{e^{\lambda_{1}},\dots, e^{\lambda_{n}}\}$, so there exist $P\in GL(n;\mathbb{R})$ such that
\begin{equation*}
  PMP^{-1}=
\begin{pmatrix}
    e^{\lambda_{1}}& & \\
    &\ddots&\\
    & &e^{\lambda_{n}}
\end{pmatrix},
\end{equation*}
and, in particular, we have 
\begin{equation}
\label{somma0} 
\sum_{i=1}^{n}\lambda_{i}=0.
\end{equation}
Consider the group homomorphism $\rho:\mathbb{C}\to GL(n;\mathbb{C})$ defined by
\begin{equation*}
   w\mapsto \begin{pmatrix}
    e^{\frac{1}{2}\lambda_{1}(w+\overline{w})}& & \\
    &\ddots&\\
    & &e^{\frac{1}{2}\lambda_{n}(w+\overline{w})}
\end{pmatrix},
\end{equation*}
and let $\mathbb{C}\ltimes_{\rho}\mathbb{C}^{n}$ be the Lie group given by the semidirect product with respect to $\rho$. Fix $\tau\in\mathbb{R}\setminus\{0\}$ and set 
\begin{equation*}
    \Gamma_{\tau}^{'}=\mathbb{Z}\oplus \sqrt{-1}\,\tau\mathbb{Z}\subset\mathbb{C},\quad \quad \Gamma^{''}_{P}=P\mathbb{Z}^{n}\oplus\sqrt{-1}P\mathbb{Z}^{n}\subset\mathbb{C}^{n},
\end{equation*}
it follows that 
\begin{equation*}
\Gamma_{P,\tau}=\Gamma_{\tau}^{'}\ltimes_{\rho}\Gamma^{''}_{P}
\end{equation*}
is a discrete subgroup of $\mathbb{C}\ltimes_{\rho}\mathbb{C}^{n}$. Therefore, the \emph{Nakamura manifold} associated to $(M,P,\tau)$ is the compact complex solvmanifold of dimension $n+1$ given by the quotient
\begin{equation*}
    N_{M,P,\tau}=\Gamma_{P,\tau}\backslash (\mathbb{C}\ltimes_{\rho}\mathbb{C}^{n}).
\end{equation*}
We will assume that  $N_{M,P,\tau}$ is not a complex torus, hence $\lambda_{i}\ne 0$ for some $i=1,\dots, n.$\\
The $(1,0)$-forms on $\mathbb{C}\ltimes_{\rho}\mathbb{C}^{n}$,
\begin{equation*}
    \varphi^{0}=d\, w,\quad\quad \varphi^{i}=e^{-\frac{1}{2}\lambda_{i}(w+\overline{w})}dz_{i}, \quad\quad i=1,\dots, n ,
\end{equation*}
descend to a global co-frame on $N_{M,P,\tau}$, for which the following structure equations hold:
\begin{equation*}
    d\varphi^{0}=0,\quad\quad d\varphi^{i}=-\frac{1}{2}\lambda_{i}(\varphi^{0}+\overline{\varphi}^{0})\wedge \varphi^{i},\quad\quad i=1,\dots, n.
\end{equation*} 
Let us introduce some notation. For a given multi-index $K=(k_1, \ldots, k_r)$, denote by $\vert K\vert$ its length, that is $\vert K\vert=r$. Let $I = (i_1, \ldots, i_p)$ and $J = (j_1, \ldots, j_q)$ be multi-indices of length $p$ and $q$ respectively, namely
\[1 \leq i_1 < \ldots < i_p \leq n, \qquad 1 \leq j_1 < \ldots < j_q \leq n,\]
and denote
\[\varphi^I = \varphi^{i_1} \wedge \ldots \wedge \varphi^{i_p}, \qquad \overline{\varphi}^J = \bar{\varphi}^{j_1} \wedge \ldots \wedge \bar{\varphi}^{j_q},\]
Set
\begin{equation}\label{cIJ}
c_{IJ} = \sum_{\sigma = 1}^p \lambda_{i_\sigma} + \sum_{\tau = 1}^q \lambda_{j_\tau},
\end{equation}
and define $f_{IJ}:\C\to\C$ as
\begin{equation}\label{fIJ}
f_{IJ}(w) = e^{-\frac{1}{2} c_{IJ} (w - \bar{w})}.   
\end{equation}
For any given pair of multi-indices $I,J$, the function $f_{IJ}$ gives rise to a well-defined complex-valued smooth function on $N_{M, P, \tau}$ if and only if
\begin{equation}\label{well-defined-fIJ}
\tau \cdot c_{IJ} \in 2\pi \cdot \Z.\\
\end{equation}
%It is straightforward to see that, if the $\del\delbar$-%lemma holds, then
%\begin{equation*}
%H_{\delbar}^{\bullet,\bullet}(N)\simeq %H_{BC}^{\bullet,\bullet}(N),
%\end{equation*}
%and, from \cite[Theorem 4.8]{CattTom}, it follows that %$N_{M, P, \tau}$ is geometrically-Bott-Chern-formal.\\
%\begin{prop}
%Let $N=N_{M,P,\tau}$ be a Nakamura manifold. If $N$ satisfies the $\del\delbar$-lemma, then $N$ is geometrically-Bott-Chern-formal.
%\end{prop}
In \cite[Theorem 4.14]{CattTom}, the validity of the $\del\delbar$-Lemma on $N_{M, P, \tau}$ is ensured by the condition:
\begin{equation}\label{eq: technical condition}
\begin{array}{c}
\text{for every $I$, $J$ it holds that}\\[5pt]
\tau \cdot c_{IJ} \in 2\pi \cdot \Z \Longleftrightarrow c_{IJ} = 0.
\end{array}
\end{equation}
Hence, let $N = N_{M, P, \tau}$ be a Nakamura manifold corresponding to a choice of $\tau$ satisfying \eqref{eq: technical condition}. Then $N$ satisfies the $\del\delbar$-Lemma and it has no K\"ahler structure.
Now we state the following.
\begin{theorem}\label{del-delbar}
Let $N = N_{M, P, \tau}$ be a Nakamura manifold satisfying \eqref{eq: technical condition}. Then, $N_{M, P, \tau}$ admits  a geometrically-Bott-Chern-formal metric.
\end{theorem}
\begin{proof}
By assumption, $N$ satisfies \eqref{eq: technical condition}; consequently, in view of \cite[Theorem 4.14]{CattTom}, $N$ satisfies the $\del\delbar$-Lemma. Therefore, 
\begin{equation*}
H_{\delbar}^{\bullet,\bullet}(N)\simeq H_{BC}^{\bullet,\bullet}(N),
\end{equation*}
and, by \cite[Theorem 4.8]{CattTom}, it follows that 
the Hermitian metric $$
g=\sum_{j=0}^n\varphi^j\otimes \overline{\varphi^j}$$
is geometrically-Bott-Chern-formal.\\ 
\end{proof}
To study obstructions to the existence of geometrically-Bott-Chern-formal metrics, we focus on a subclass of Nakamura manifolds for which the $\del\delbar$-Lemma does not hold. We make the following assumption:
\begin{equation}\label{8mod}
\begin{array}{c}
\text{ there exist multi-indices $I$, $J$ with $I\cap J= \emptyset$ such that}\\[5pt]
\tau \cdot c_{IJ} \in 2\pi \cdot \Z\setminus\{0\},
\end{array}
\end{equation}
which leads to the following result:
\begin{theorem}\label{gNakamura}
Let $N=N_{M,P,\tau}$ be a Nakamura manifold corresponding to a choice of $\tau$ satisfying \eqref{8mod}. Then $N$ admits no geometrically-Bott-Chern-formal metrics.
\end{theorem}
Before giving the proof of Theorem \ref{gNakamura} we will state and prove Lemma \ref{BC-coh} and Lemma \ref{powers}. In particular, Lemma \ref{BC-coh} provides an explicit description of the Bott--Chern cohomology generators on $N$.
%The following Lemma \ref{BC-coh} and Lemma \ref{powers} %will be useful for the proof of Theorem \ref{gNakamura}. %In particular, Lemma \ref{BC-coh} provides an explicit %description of the Bott--Chern cohomology generators on %$N$.
\begin{lemma}\label{BC-coh}
Let $N=N_{M,P,\tau}$ be a Nakamura manifold. Then the $(p,q)$-Bott-Chern cohomology $H_{BC}^{p,q}(N)$ is isomorphic to the complex vector space generated by\\
$$
\begin{array}{llll}
    \delta_{c_{LM}}f_{LM}\varphi^{L\overline{M}},&  \delta_{c_{LM}}\overline{f_{LM}}\varphi^{L\overline{M}}& \quad\quad  |L|=p,& |M|=q\\[5pt]
    f_{LM}\varphi^{0L\overline{M}},& \delta_{c_{LM}} \overline{f_{LM}}\varphi^{0L\overline{M}}& \quad\quad  |L|=p-1, &|M|=q \\[5pt]
     \delta_{c_{LM}}f_{LM}\varphi^{\overline{0}L\overline{M}},& \overline{f_{LM}}\varphi^{\overline{0}L\overline{M}}& \quad\quad  |L|=p,& |M|=q-1\\[5pt]
f_{LM}\varphi^{0\overline{0}L\overline{M}},&  \overline{f_{LM}}\varphi^{0\overline{0}L\overline{M}} & \quad\quad  |L|=p-1,& |M|=q-1
\end{array}
$$
where 
\begin{equation*}
 \delta_{c_{LM}}=
\begin{cases}
 1& \text{if}\quad c_{LM}=0,\\
 0&\text{if}\quad c_{LM}\ne 0,
\end{cases}
\end{equation*}
and $L,\, M$ correspond to multi-indices such that $\tau \cdot c_{LM}\in 2\pi\mathbb{Z}$. Each generator is Bott-Chern harmonic with respect to the Hermitian metric
$$
g=\sum_{j=0}^n\varphi^j\otimes \overline{\varphi^j}
$$
on $N$.
\end{lemma}
\begin{proof}[Proof of Lemma \ref{BC-coh}]
  By \cite[Theorem 2.16]{AngKas}, we have
\begin{equation*}
    H_{BC}^{p,q}(N)\simeq H_{BC}^{p,q}(C^{p,q})
\end{equation*}
where $C^{p,q}=B^{p,q}+\overline{B}^{p,q}$, and by \cite[Lemma 4.6]{CattTom} $B^{p,q}$ is generated by
$$
\begin{array}{lll}
    f_{LM}\varphi^{L\overline{M}}& \quad\quad  |L|=p,& |M|=q\\[5pt]
    f_{LM}\varphi^{0L\overline{M}}& \quad\quad  |L|=p-1, &|M|=q \\[5pt]
    f_{LM}\varphi^{\overline{0}L\overline{M}}& \quad\quad  |L|=p,& |M|=q-1\\[5pt]
f_{LM}\varphi^{0\overline{0}L\overline{M}} & \quad\quad  |L|=p-1,& |M|=q-1
\end{array}
$$
corresponding to the multi-indices $L, M$ such that the functions $f_{LM}$ are well-defined on $N$, that is, if and only if $\tau \cdot c_{LM} \in 2\pi \mathbb{Z}$.\\
Since $\delbar B^{p,q}=0$, it follows 
\begin{equation*}
     H_{BC}^{p,q}(N)\simeq \text{Ker}(d: C^{p,q}\to C^{p+1, q+1})\cap C^{p,q}
\end{equation*} 
Finally, a straightforward computation yields the set of Bott-Chern harmonic generators stated in the lemma.
\end{proof}
\begin{rem}\label{rem-1}
As a consequence, one can easily compute the Aeppli cohomology by taking the Hodge star operator of the Bott-Chern harmonic generators listed in Lemma \ref{BC-coh}. For the sake of completeness, we explicitly present them below:
$$
\begin{array}{llll}
    f_{LM}\varphi^{L\overline{M}},&  \overline{f_{LM}}\varphi^{L\overline{M}}& \quad\quad  |L|=p,& |M|=q\\[5pt]
    \delta_{c_{LM}}  f_{LM}\varphi^{0L\overline{M}},&\overline{f_{LM}}\varphi^{0L\overline{M}}& \quad\quad  |L|=p-1, &|M|=q \\[5pt]
    f_{LM}\varphi^{\overline{0}L\overline{M}},&  \delta_{c_{LM}}\overline{f_{LM}}\varphi^{\overline{0}L\overline{M}}& \quad\quad  |L|=p,& |M|=q-1\\[5pt]
\delta_{c_{LM}}f_{LM}\varphi^{0\overline{0}L\overline{M}},&  \delta_{c_{LM}}\overline{f_{LM}}\varphi^{0\overline{0}L\overline{M}} & \quad\quad  |L|=p-1,& |M|=q-1.
\end{array}
$$
corresponding to the multi-indices $L, M$ for which the functions $f_{LM}$ are well-defined on $N$.
\end{rem}
Finally, the following lemma shows that the integrals of all non-trivial powers of the functions $f_{IJ}$, defined in \eqref{fIJ}, vanish on $N$. More precisely,  
\begin{lemma}\label{powers}
Let $N = N_{M, P, \tau}$ be a Nakamura manifold, and let $I,J$ be multi-indices such that \[\tau\cdot c_{IJ} \in 2\pi\mathbb{Z}\setminus\{0\}.\] Then, for all $k\in \mathbb{Z}\setminus\{0\}$, we have 
\begin{equation*}
   \int_{N}f^{k}_{IJ}\,\Omega=0,
\end{equation*}
where $\Omega=\frac{\omega^{n+1}}{(n+1)!}$ is the standard volume form and $\omega$ is the fundamental form of the Hermitian metric 
\[g=\sum_{j=0}^n\varphi^j\otimes \overline{\varphi^j}.
\]
\end{lemma}
\begin{proof}[Proof of Lemma \ref{powers}]
A direct computation gives
    \begin{align*}
    &\int_{N}f^{k}_{IJ}\,\Omega\\
    &=|\det P|^{2} \int_{0}^{\tau}e^{-k\sqrt{-1}c_{IJ}y}\,dy\\
    &=|\det P|^{2}\int_{0}^{\tau}\big(\cos(kc_{IJ}y)-\sqrt{-1}\sin(kc_{IJ}y)\big)\,dy\\
    &=0,
 \end{align*}
 where we set $w=x+\sqrt{-1}y\in\mathbb{C}$, and we denote by $|\det P|^2$ the square modulus of the determinant of the matrix $P$. The integral vanishes since $\tau \cdot c_{IJ} \in 2\pi \mathbb{Z}\setminus\{0\}$.
\end{proof}
We are now ready to prove Theorem~\ref{gNakamura}.
\begin{proof}[Proof of Theorem \ref{gNakamura}]
We will prove that $N$ does not admit any geometrically-Bott-Chern-formal metric by exhibiting a non-vanishing triple $ABC$-Massey product on $N$. \\By assumption \eqref{8mod}, there exist multi-indices $I, J$ of length $|I|=p$, $|J|=q$, with $I\cap J=\emptyset$, such that
\begin{equation*}
    \tau\cdot c_{IJ}\in 2\pi\mathbb{Z}\setminus\{0\}.
\end{equation*}
For such a pair $(I,J)$, the function $f_{IJ}$ is a non-constant complex-valued smooth function on $N$.\newline
Then, let us consider the following classes in Bott-Chern cohomology, defined respectively as
$$
\mathfrak{a}_{12}=[\alpha_{12}]=[f_{IJ}\,\varphi^{0IJ}]\in H_{BC}^{p+q+1, 0}(N),\,\mathfrak{a}_{23}=[\alpha_{23}]=[c_{IJ}^{2}\overline{f_{IJ}}\,\varphi^{\overline{0IJ}}]\in H_{BC}^{0, p+q+1}(N).
$$
Setting $T=(I\cup J)^{c}$, where the symbol $c$ indicates the complement of a set, consider the form $f_{IJ}\,\varphi^{\overline{0T}}$. By the definition of $c_{IJ}$ (see \eqref{cIJ}), since $\sum_{j=1}^{n}\lambda_j=0$ and $I\cap\, J = \emptyset$, it follows that
$$
c_{T}=-c_{IJ},
$$
hence a direct computation shows that it defines a Bott--Chern cohomology class. We can then set
$$
\mathfrak{a}_{34}=[\alpha_{34}]=[f_{IJ}\,\varphi^{\overline{0T}}]\in H_{BC}^{0, n-p-q+1}(N).
$$
Since 
$$
\alpha_{12}\wedge \alpha_{23}=\del\delbar\,(-1)^{p+q}\varphi^{IJ\overline{IJ}}\,,
\quad \alpha_{23}\wedge \alpha_{34}=0,
$$
we have
$$\mathfrak{a}_{12}\cup \mathfrak{a}_{23}=0,\quad \mathfrak{a}_{23}\cup\mathfrak{a}_{34}=0.$$
Therefore, the following $ABC$-Massey product is well-defined:
\begin{eqnarray*}
\langle\mathfrak {a}_{12},\mathfrak {a}_{23},\mathfrak {a}_{34}\rangle_{ABC}&=&[(-1)^{p+q+1}\varphi^{IJ\overline{IJ}}\wedge f_{IJ}\,\varphi^{\overline{0T}}]=[-f_{IJ}\,\varphi^{IJ\overline{0IJT}}]\\[5pt]
&\in& \frac{H_A^{p+q,n+1}(N)}{H_A^{p+q, p+q}(N)\cup[f_{IJ}\,\varphi^{\overline{0T}}]}.
\end{eqnarray*}
We will now show that 
$$
[-f_{IJ}\,\varphi^{IJ\overline{0IJT}}]\notin H_A^{p+q, p+q}(N)\cup[f_{IJ}\,\varphi^{\overline{0T}}],
$$
so that the triple $ABC$-Massey product $\langle\mathfrak {a}_{12},\mathfrak {a}_{23},\mathfrak {a}_{34}\rangle_{ABC}$ is non-vanishing on $N$.\\\\
By contradiction, suppose there exists an Aeppli harmonic form
$$\alpha\in \mathcal{A}^{p+q,p+q}(N),$$
and forms $R\in \mathcal{A}^{p+q-1,n+1}(N),\quad S\in \mathcal{A}^{p+q,n}(N)$ such that\\
\begin{equation}\label{eqfij}
-f_{IJ}\,\varphi^{IJ\overline{0IJT}}=\alpha\wedge f_{IJ}\, \varphi^{\overline{0T}} +\del R^{p+q-1,n+1}+\delbar S^{p+q,n}.
\end{equation}
\\
Let $g=\sum_{j=0}^n\varphi^j\otimes \overline{\varphi^j}$ be the Hermitian metric on $N$, we denote by $\omega$ its associated fundamental form and by $\Omega=\frac{\omega^{n+1}}{(n+1)!}$ the standard volume form.\\
Multiplying equation \eqref{eqfij} by $\ast(-f_{IJ}\,\varphi^{IJ\overline{0IJT}})$, we get\\ 
\begin{equation}\label{eqfij2}
    C\,\Omega= \ast(-f_{IJ}\,\varphi^{IJ\overline{0IJT}})\wedge\big( \alpha\wedge f_{IJ}\, \varphi^{\overline{0T}} +\del R^{p+q-1,n+1}+\delbar S^{p+q,n} \big),
\end{equation}
\\where $C=g(f_{IJ}\,\varphi^{IJ\overline{0IJT}}, \,f_{IJ}\,\varphi^{IJ\overline{0IJT}})\in\mathbb{C}\setminus\{0\}.$\\
By Remark \ref{rem-1}, we can write
\begin{equation*}
\ast(-f_{IJ}\,\varphi^{IJ\overline{0IJT}})\wedge\big( \alpha\wedge f_{IJ}\, \varphi^{\overline{0T}})=(\mu_{1}f_{IJIJ}+\mu_{2}\overline{f_{IJIJ}})\,\Omega=(\mu_{1}f^{2}_{IJ}+\mu_{2}\overline{f^{2}_{IJ}})\,\Omega
\end{equation*}
for some $\mu_{1},\,\mu_{2}\in\mathbb{C}$.\\
Therefore, equation \eqref{eqfij2} can be rewritten as
\begin{equation*}
    C\,\Omega= (\mu_{1}f^{2}_{IJ}+\mu_{2}\overline{f^{2}_{IJ}})\,\Omega+\ast(-f_{IJ}\,\varphi^{IJ\overline{0IJT}})\wedge \big(\del R^{p+q-1,n+1}+\delbar S^{p+q,n} \big).
\end{equation*}
Since 
%T=(I\cup J)^c$, by the definition of $c_{IJ}$ see \eqref{cIJ}, taking into account that $\sum_{j=1}^{n}\lambda_j=0$, it turns out that
$$
c_{T}=-c_{IJ},
$$
a direct computation shows $\ast(f_{IJ}\,\varphi^{IJ\overline{0IJT}})$ is Aeppli-harmonic, which implies 
$$d^{\ast}(f_{IJ}\,\varphi^{IJ\overline{0IJT}})=0,$$ as a result, we can rewrite the last equation as
\begin{equation*}
    C\,\Omega= (\mu_{1}f^{2}_{IJ}+\mu_{2}\overline{f^{2}_{IJ}})\,\Omega+d\, P
\end{equation*}
where 
$$P=\ast(-f_{IJ}\,\varphi^{IJ\overline{0IJT}})\wedge (-1)^{n-p-q+1}\,\big( R^{p+q-1,n+1}+ S^{p+q,n} \big).$$
Integrating on $N$ and applying Stokes' Theorem, it becomes
\begin{equation*}
    C\,\int_{N}\Omega= \mu_{1}\,\int_{N} f^{2}_{IJ}\,\Omega+\mu_{2}\,\int_{N}\overline{f^{2}_{IJ}}\,\Omega,
\end{equation*}
and by applying Lemma \ref{powers} we obtain 
\begin{equation*}
     C\,\int_{N}\Omega=0.
\end{equation*}
This is a contradiction, as $C$ is a non-zero constant and $\Omega$ is a volume form.\\

\end{proof}

Concerning the existence of astheno-K\"{a}hler metrics, we prove the following result: 
\begin{prop}\label{no-astheno}
Let $N=N_{M,P,\tau}$ be a Nakamura manifold. Then $N$ does not admit any astheno-K\"{a}hler metric. 
\end{prop}
\begin{proof}
Consider  the $(2,2)$-form defined by $$\phi^{2,2}:=\varphi^{0i\overline{0i}},$$ for some $1\le i\le n$ such that $\lambda_{i}\ne 0$. Such an index exists because $N$ is not a complex torus. A direct computation shows that
\begin{equation*}
\phi^{2,2}=\,-\frac{1}{\lambda_{i}^{2}}\del\delbar\varphi^{i\overline{i}},
\end{equation*}
and hence the thesis follows from \cite[Lemma 3.5]{STa}. For the sake of completeness, we recall here the argument: suppose that $\omega$ is the fundamental form of an astheno-K\"ahler metric on $N$, so that $$\del\delbar\omega^{n-1}=0,$$
then we would obtain
\begin{eqnarray*}
0&< & \int_{N}\omega^{n-1}\wedge \phi^{2,2}= -\frac{1}{\lambda_{i}^{2}}\int_{N}\omega^{n-1}\wedge \del\delbar\,\varphi^{i\overline{i}}= -\frac{1}{\lambda_{i}^{2}} \int_{N}\del\delbar\,\omega^{n-1}\wedge\varphi^{i\overline{i}}\,=\,0,
\end{eqnarray*}
which is a contradiction.
\end{proof}
\section{Non-formality of Compact Quotients of  $\mathbb{C}^{2n}\ltimes_{\rho} \mathbb{C}^{2m}$}\label{LT}
We recall the construction of a family of compact complex solvmanifolds, not of completely solvable type, obtained as compact quotients of a solvable Lie group of the form $\mathbb{C}^{2n}\ltimes_{\rho} \mathbb{C}^{2m}$, as provided in \cite{LuTom}.\\
Let $M\in SL(2m;\mathbb{Z})$ be a diagonalizable matrix over $\mathbb{R}$ such that
\begin{equation*}
  PMP^{-1}=
\begin{pmatrix}
    e^{\lambda}& & & & & \\
    &e^{-\lambda}& & & & \\
    & & \ddots& & &\\
      & && &e^{\lambda} &\\
   & & & & & e^{-\lambda}
\end{pmatrix},
\end{equation*}
for some $\lambda\in \mathbb{R}\setminus\{0\}$ and $P\in GL(2m;\mathbb{R})$.\\
Let $w=(w_{1}, \dots, w_{2n})$ denote the coordinates of $\mathbb{C}^{2n}$, and consider the group homomorphism $$\rho:\mathbb{C}^{2n}\to GL(2m;\mathbb{C}),$$ defined by
\begin{equation*}
 w\mapsto \begin{pmatrix}
    e^{\phi(w)}&  &&   \\
    &e^{-\overline{\phi(w)}}&   & \\
    & \ddots&  &\\
      &  &e^{\phi(w)} &\\
   &  & & e^{-\overline{\phi(w)}}
\end{pmatrix},
\end{equation*}
where $$\phi(w)=\lambda(w_{1}+\overline{w}_{2}+\dots +w_{2n-1}+\overline{w}_{2n}),$$
so we can consider the Lie group $\mathbb{C}^{2n}\ltimes_{\rho}\mathbb{C}^{2m}$ given by the semidirect product with respect to $\rho$.\\
Fix the non-zero diagonal matrix $\underline{\mu}= \text{diag}(\mu_{1},\dots ,\mu_{2n})$, and define
\begin{equation*}
    \Gamma_{\underline{\mu}}^{'}=\mathbb{Z}^{2n}\oplus \sqrt{-1}\,\underline{\mu}\,\mathbb{Z}^{2n},\quad \quad \Gamma^{''}_{P}=P\,\mathbb{Z}^{2m}\oplus\sqrt{-1}P\,\mathbb{Z}^{2m}.
\end{equation*}
If $\underline{\mu}=\text{diag}(\mu_{1},\dots,\mu_{2n})$ has entries of the form $\mu_{i}=\frac{2k_{i}\pi}{\lambda}$ with $k_{i}\in\mathbb{Z}$, then 
\begin{equation*}
\Gamma_{\underline{\mu}, P}=\Gamma_{\underline{\mu}}^{'}\ltimes_{\rho}\Gamma^{''}_{P}
\end{equation*}
is a discrete subgroup of $\mathbb{C}^{2n}\ltimes_{\rho}\mathbb{C}^{2m}$. Hence, associated with the choice of $\underline{\mu}$ and $P$, we can consider the compact complex solvmanifold of dimension $2n+2m$ given by the quotient
\begin{equation*}
    N_{\underline{\mu}, P}=\Gamma_{\underline{\mu}, P}\backslash (\mathbb{C}^{2n}\ltimes_{\rho}\mathbb{C}^{2m}).
\end{equation*}
The $(1,0)$-forms on $\mathbb{C}^{2n}\ltimes_{\rho}\mathbb{C}^{2m}$:
\begin{equation*}
    \varphi^{i}=d\, w_{i}\qquad i=1,\dots, 2n,
\end{equation*}
and 
\begin{equation*}
 \psi^{2j+1}=e^{-\lambda(w_{1}+\overline{w}_{2}+\dots+w_{2n-1}+\overline{w}_{2n})}dz_{2j+1}, \,\,  \psi^{2j+2}=e^{\lambda(\overline{w}_{1}+w_{2}+\dots+\overline{w}_{2n-1}+w_{2n})}dz_{2j+2},
\end{equation*}
for $j=0,\dots, m-1$, define a global co-frame on $N=N_{\underline{\mu},P}$. Setting $$\eta=\varphi^{1}+\overline{\varphi}^{2}+\dots +\varphi^{2n-1}+\overline{\varphi}^{2n},$$ the structure equations are given by:
\begin{equation*}
    d\varphi^{i}=0,\quad\quad d\psi^{2j+1}=-\lambda\,\eta\wedge \psi^{2j+1}\quad\quad d\psi^{2j+2}=\lambda\,\overline{\eta}\wedge \psi^{2j+2}.
\end{equation*} We will prove the following:
\begin{theorem}\label{Lu1}
Let $N=N_{\underline{\mu}, P}$. Then $N$ is not geometrically-Bott-Chern-formal.
\end{theorem}
\begin{proof}
According to \cite[Theorem 2.4]{AT-1}, in order to show that $N$ does not admit any geometrically-Bott-Chern-formal metric, we will construct a non-vanishing triple $ABC$-Massey product.\\
Let
\begin{equation*}
    \sigma=\varphi^{1}+\varphi^{2}+\dots+\varphi^{2n},
\end{equation*}
and define
\begin{equation*}
\beta_{j}=
\begin{cases}
  e^{-\lambda(w_{2}-\overline{w}_{2}+\dots+ w_{2n}-\overline{w}_{2n})}\qquad\quad\,  j \,\,\text{odd},\\
  e^{\lambda(w_{1}-\overline{w}_{1}+\dots+ w_{2n-1}-\overline{w}_{2n-1})}\qquad  j\,\, \text{even}.
\end{cases}
\end{equation*}
Note that the complex-valued function $\beta_{j}$ is well-defined on $N$, since it is trivial on $\Gamma_{\underline{\mu}, P}$ by the choice of $\underline{\mu}$. We will use the notation $\varphi^{A\overline{B}}=\varphi^{A}\wedge\overline{\varphi}^{B}$.\\\\
Consider the following Bott-Chern cohomology classes:
\begin{align*}
\mathfrak{a}_{12}=[\alpha_{12}]=[\beta_{2}\,&\sigma\wedge\psi^{2}]\in H_{BC}^{2,0}(N), \quad\mathfrak{a}_{23}=[\alpha_{23}]=[-\lambda^{2}\overline{\beta}_{1}\,\overline{\sigma}\wedge\psi^{\overline{1}}]\in H_{BC}^{0,2}(N),\\
&\mathfrak{a}_{34}=[\alpha_{34}]=[\overline{\beta}_{2}\,\overline{\sigma}\wedge\psi^{\overline{2}}]\in H_{BC}^{0,2}(N).
\end{align*}
Indeed, for example, by a direct computation one can verify that 
$$
\begin{array}{lll}
    d\alpha_{12}&=&\\[5pt]
    &=&\lambda\,\beta_{2}(\varphi^{1}-\varphi^{\overline{1}}+\dots+\varphi^{2n-1}-\varphi^{\overline{2n-1}})\wedge\sigma\wedge\psi^{2}\\[5pt]
    &+&\lambda\,\beta_{2}(\varphi^{\overline{1}}
+\varphi^{2}+\dots+\varphi^{\overline{2n-1}}+\varphi^{2n})\wedge\sigma\wedge\psi^{2}\\[5pt]
&=&\lambda\,\beta_{2}\,\sigma\wedge\sigma\wedge\psi^{2}=0,
\end{array}
$$
and thus $\mathfrak{a}_{12}\in H_{BC}^{2,0}(N)$.\\ 
Now, observe that:
$$
\alpha_{12}\wedge \alpha_{23}=\del\delbar\,\,\overline{\beta}_{1}\,\beta_{2}\,\psi^{\overline{1}2}\,
\quad \alpha_{23}\wedge \alpha_{34}=0,
$$
and thus,
$$\mathfrak{a}_{12}\cup \mathfrak{a}_{23}=0,\quad \mathfrak{a}_{23}\cup\mathfrak{a}_{34}=0.$$
Therefore, the following $ABC$-Massey product is well-defined:
\begin{eqnarray*}
\langle\mathfrak {a}_{12},\mathfrak {a}_{23},\mathfrak {a}_{34}\rangle_{ABC}&=&[-\overline{\beta}_{1}\,\beta_{2}\,\psi^{\overline{1}2}\,\overline{\beta}_{2}\,\overline{\sigma}\wedge\psi^{\overline{2}}]=[\overline{\beta}_{1}\,\overline{\sigma}\wedge\psi^{2\overline{12}}]\\[5pt]
&\in& \frac{H_A^{1,3}(N)}{H_{A}^{1,1}(N)\cup[\overline{\beta}_{2}\,\overline{\sigma}\wedge\psi^{\overline{2}}]}.
\end{eqnarray*}
We now show that 
$$
[\overline{\beta}_{1}\,\overline{\sigma}\wedge\psi^{2\overline{12}}]\not\in H_{A}^{1,1}(N)\cup[\overline{\beta}_{2}\,\overline{\sigma}\wedge\psi^{\overline{2}}],
$$
which implies that the triple $ABC$-Massey product $\langle\mathfrak {a}_{12},\mathfrak {a}_{23},\mathfrak {a}_{34}\rangle_{ABC}$ is non-vanishing on $N$.\\
Suppose, by contradiction, that there exist $\xi^{1,1}\in\mathcal{A}^{1,1}(N)$ with 
$\del\delbar\xi^{1,1}=0$,
 $R^{0,3}\in\mathcal{A}^{0,3}(N)$ and $S^{1,2}\in\mathcal{A}^{1,2}(N),$ such that
\begin{align}\label{eqq1}
\overline{\beta}_{1}\,\overline{\sigma}\wedge\psi^{2\overline{12}}=\overline{\beta}_{2}\,\xi^{1,1}\wedge\overline{\sigma}\wedge\psi^{\overline{2}}+\del R^{0,3}+\delbar S^{1,2}. 
\end{align}
Let $g$ be the invariant Hermitian metric on $N$ whose associated fundamental form is given by 
$$\omega=\frac{\sqrt{-1}}{2}\sum_{j=1}^{2n} \varphi^{j}\wedge\varphi^{\overline{j}}+\frac{\sqrt{-1}}{2}\sum_{j=1}^{2m} \psi^{j}\wedge\psi^{\overline{j}},$$ and let $\Omega=\frac{\omega^{2n+2m}}{(2n+2m)!}$ denote the standard volume form.\\  Multiplying equation \eqref{eqq1} by $\ast(\overline{\beta}_{1}\,\overline{\sigma}\wedge\psi^{2\overline{12}})$, we get \\
\begin{equation}\label{bbb}
    C\,\Omega= \ast(\overline{\beta}_{1}\,\overline{\sigma}\wedge\psi^{2\overline{12}})\wedge (\overline{\beta}_{2} \,\xi^{1,1}\wedge\overline{\sigma}\wedge\psi^{\overline{2}}+\del R^{0,3}+\delbar S^{1,2} ),
\end{equation}
\\where $C=g(\overline{\beta}_{1}\,\overline{\sigma}\wedge\psi^{2\overline{12}},\,\overline{\beta}_{1}\,\overline{\sigma}\wedge\psi^{2\overline{12}})\in\mathbb{C}\setminus\{0\}.$\\
Observing that
\begin{equation*}
\overline{\sigma}\wedge\sum_{i=1}^{2n}\varphi^{\overline{1}\cdots\overline{\hat{i}}\cdots\overline{2n}}=0,
\end{equation*}
where the hat indicates the omission of the corresponding index, we obtain 
\begin{equation*}
\ast(\overline{\beta}_{1}\,\overline{\sigma}\wedge\psi^{2\overline{12}})\wedge (\overline{\beta}_{2}\,\xi^{1,1}\wedge\overline{\sigma}\wedge\psi^{\overline{2}})
=0.
\end{equation*}
Therefore, equation \eqref{bbb} can be rewritten as
\begin{equation*}
    C\,\Omega= \ast(\overline{\beta}_{1}\,\overline{\sigma}\wedge\psi^{2\overline{12}})\wedge (\del R^{0,3}+\delbar S^{1,2} ). 
\end{equation*}
 A direct computation shows that $\overline{\beta}_{1}\,\overline{\sigma}\wedge\psi^{2\overline{12}} $ is Aeppli-harmonic, which implies 
$$d^{\ast}(\overline{\beta}_{1}\,\overline{\sigma}\wedge\psi^{2\overline{12}})=0,$$ as a result, we can rewrite the last equation as
\begin{equation*}
    C\,\Omega=d\, P
\end{equation*}
where $P=\ast(\overline{\beta}_{1}\,\overline{\sigma}\wedge\psi^{2\overline{12}})\wedge( R^{0,3}+ S^{1,2})$.\\
Integrating on $N$ and applying Stokes' Theorem, 
we obtain a contradiction, as $C$ is a non-zero constant and  $\Omega$ is a volume form.\\
Hence, we have shown that 
\begin{equation*}
\langle\mathfrak {a}_{12},\mathfrak {a}_{23},\mathfrak {a}_{34}\rangle_{ABC}=[\overline{\beta}_{1}\,\overline{\sigma}\wedge\psi^{2\overline{12}}]\ne 0
\end{equation*}
which means $N$ admits a non-vanishing $ABC$-Massey product.\\
\end{proof}
On astheno-K\"{a}hler metrics, we have the following:
\begin{prop}\label{Lu2}
Let $N=N_{\underline{\mu},P}$. Then $N$ does not admit any astheno-K\"{a}hler metric. 
\end{prop}
\begin{proof}
Suppose, by contradiction, that there exists an astheno-K\"{a}hler metric on $N$, then the associated fundamental form $\omega$ satisfies 
\begin{equation*}
    \del\delbar\omega^{2n+2m-2}=0.
\end{equation*}
Let
\begin{equation*}
    \sigma=\varphi^{1}+\varphi^{2}+\dots+\varphi^{2n},
\end{equation*}
and consider the $\del\delbar$-exact $(2,2)$-form
\begin{equation*}
\phi^{2,2}:=\sigma\wedge\psi^{1}\wedge\overline{\sigma}\wedge\psi^{\overline{1}}=-\frac{1}{\lambda^{2}}\del\delbar\, \psi^{1\overline{1}}.
\end{equation*}
%where the complex-valued function $e^{\lambda(w_{1}-\overline{w}_{1}+w_{2}-\overline{w}_{2})}$ is well-defined, since it is trivial on $\Gamma_{\underline{\mu}, P}$ by the choice of %$\underline{\mu}$. 
%Observe that 
%\begin{equation*}
 %   \omega^{2n+2m-2}\wedge\phi^{2,2}=\,e^{\lambda(w_{1}-\overline{w}_{1}+w_{2}-\overline{w}_{2})}f\, \Omega
%\end{equation*}
%where $f$ is a positive real smooth function on $N$ and $\Omega=\frac{\omega^{2n+2m}}{(2n+2m)!}$ denote the standard volume form on $N$. 
It follows
\begin{align*} 
0<\int_{N}&\omega^{2n+2m-2}\wedge\phi^{2,2}\\
&=-\frac{1}{\lambda^{2}}\int_{N}\omega^{2n+2m-2}\wedge\del\delbar \psi^{1\overline{1}}\\
&=-\frac{1}{\lambda^{2}}\int_{N}\del\delbar(\omega^{2n+2m-2})\wedge\psi^{1\overline{1}}=0.
\end{align*}
This is absurd.
\end{proof}


\begin{thebibliography}{12}

%\bibitem{Al1} L. Alessandrini, Proper Modifications of generalized %$p$-K\"ahler manifolds, {\em J. Geom. Anal.} {\bf 27} (2017), 947-%-967,  https://doi.org/10.1007/s12220-016-9705-z 
\bibitem{AngKas}  D. Angella, H. Kasuya, Bott–Chern cohomology of solvmanifolds, \emph{Ann. Glob. Anal. Geom.} {\bf 52} (2017), 363--411, https://doi.org/10.1007/s10455-017-9560-6
\bibitem{AngKas2} D. Angella, H. Kasuya, Cohomologies of deformations of solvmanifolds and closedness of some properties, {\em North-West. Eur. J. Math} {\bf 3} (2017), 75--105, {\tt arXiv:1305.6709 math[CV]}.
\bibitem{AN} R.M. Arroyo, M. Nicolini, SKT structures on nilmanifolds, \emph{Math. Z.} {\bf 302} (2022), 1307--1320, https://doi.org/10.1007/s00209-022-03107-3
\bibitem{ASfe} D. Angella, T. Sferruzza, Geometric Formalities Along the Chern-Ricci Flow, {\em Complex Anal. Oper. Theory} {\bf 14} (2020), https://doi.org/10.1007/s11785-019-00971-6
\bibitem{AT-1} D. Angella, A. Tomassini, On Bott–Chern cohomology and formality, \emph{J. Geom. Phys.} \textbf{93} (2015), 52--61, https://doi.org/10.1016/j.geomphys.2015.03.004
\bibitem{trivial-bundle}  M.L. Barberis, I. Dotti, M. Verbitsky, Canonical bundles of complex nilmanifolds, with applications to hypercomplex geometry, \emph{Math. Res. Lett.} \textbf{16} no. 2 (2009), 331-347, https://dx.doi.org/10.4310/MRL.2009.v16.n2.a10
\bibitem{BigRol} L. Bigalke, S. Rollenske, Erratum to: The Fr\"olicher spectral sequence can be arbitrarily non-degenerate,
{\em Math. Ann.} {\bf 358} (2014), 1119--1123, https://doi.org/10.1007/s00208-013-0996-0
 \bibitem{CattTom} A. Cattaneo, A. Tomassini, $\del\delbar$-lemma and $p$-K\"{a}hler structures on families of solvmanifolds, {\em Math. Z.} {\bf 308} no. 56 (2024), https://doi.org/10.1007/s00209-024-03612-7
 \bibitem{JNilp} L. A. Cordero, M. Fern\'andez, A. Gray, L. Ugarte, Compact nilmanifolds with nilpotent complex structures: Dolbeault
cohomology, \emph{Trans. Amer. Math. Soc.} {\bf 352} no. 12 (2000), 5405–5433, http://www.jstor.org/stable/221897
\bibitem{FV} A. Fino, L. Vezzoni, On the existence of balanced and SKT metrics on nilmanifolds, \emph{Proc. Amer. Math. Soci.} {\bf 144} (2016), 2455--2459,
http://dx.doi.org/10.1090/proc/12954
\bibitem{HK} R. Harvey, A. W. Knapp, Positive $(p,p)$ forms, Wirtinger's inequality, and currents, in \emph{Value distribution theory (Proc. Tulane Univ. Program, Tulane Univ., New Orleans, La., 1972-1973)}, Part A, vol. {\bf 25} of \emph{Pure Appl. Math.}, Dekker, New York, 1974, 43--62.
\bibitem{Has05} K. Hasegawa, Complex and K\"ahler structures on compact solvmanifolds, {\em J. Symplect. Geom.} {\bf 3} no. 4 (2005), 749--767.
\bibitem{Has06} K. Hasegawa, A note on compact solvmanifolds with K\"ahler structures, \emph{Osaka J. Math.} {\bf 43} no. 1 (2006), 131--135.
\bibitem{Kas} H. Kasuya, Techniques of computations of Dolbeault cohomology
of solvmanifolds, \emph{Math. Z.} {\bf 273} (2013), 437--447, https://doi.org/10.1007/s00209-012-1013-0
%\bibitem{hirzebruch} F. Hirzebruch, Some problems on differentiable and complex manifolds, \emph{Ann. Math. (2)} \textbf{60}, (1954). 213--236.

\bibitem{Kot} D. Kotschick, On products of harmonic forms, \emph{Duke Math. J.} {\bf 107} no. 3 (2001), 521--531, https://doi.org/10.1215/s0012-7094-01-10734-5
\bibitem{LuTom} F. Lusetti, A. Tomassini, Hard Lefschetz Condition on symplectic non-K\"{a}hler solvmanifolds, {\em Math. Z.} {\bf 311} no. 74 (2025), https://doi.org/10.1007/s00209-025-03878-5
%\bibitem{Massey} W.S. Massey, Some higher order cohomology %operations, {\em Sympos. Int. Topologia Algebraica} (1958), 145-154
\bibitem{MS} A. Milivojević, J. Stelzig, Bigraded notions of formality and Aeppli--Bott--Chern--Massey products, \emph{Comm. Anal. Geom.} {\bf 32} (2024), 2901-2933, 
https://doi.org/10.4310/CAG.241231035705
\bibitem{Nak} I. Nakamura, Complex parallelisable manifolds and their small deformations, {\em J. Differ. Geom.} {\bf 10} no. 1 (1975),
85--112, https://doi.org/10.4310/jdg/1214432677
%\bibitem{Neis} J. Neisendorfer, L. Taylor, Dolbeault homotopy %theory, {\em Trans. Amer. Math. Soc.} {\bf 245} (1978), 183-210, %https://doi.org/10.1090/S0002-9947-1978-0511405-5
\bibitem{PSZ24} G. Placini, J. Stelzig, L. Zoller, Nontrivial Massey products on compact K\"ahler manifolds, {\tt arXiv:2404.09867 [math.AT]} (2024).
\bibitem{R} L. Rubini, Some computations on trivial canonical-bundle solvmanifolds, {\em J.Geom. Phys.} {\bf 216} (2025), 105586, https://doi.org/10.1016/j.geomphys.2025.105586
\bibitem{Schw} M. Schweitzer, Autour de la cohomologie de Bott-Chern, $\mathtt{arXiv:0709.3528v1\, [math.AG]}$ (2007).
\bibitem{SfeTar} T. Sferruzza, N. Tardini, $p$-K\"ahler and balanced structures on nilmanifolds with nilpotent complex structures, {\em Ann. Glob. Anal. Geom.} {\bf 62} (2022), 869--881, https://doi.org/10.1007/s10455-022-09867-9
\bibitem{SfeTom1} T. Sferruzza, A. Tomassini, Dolbeault and Bott-Chern formalities: deformations and $\del\delbar$-lemma, \emph{J. Geom. Phys.} {\bf 175} (2022), 104470, https://doi.org/10.1016/j.geomphys.2022.104470
\bibitem{STa} T. Sferruzza, A. Tomassini, On cohomological and formal properties of strong K\"ahler with torsion and astheno-K\"ahler metrics, {\em Math. Z.} {\bf 304} no. 55 (2023), https://doi.org/10.1007/s00209-023-03303-9
\bibitem{STb} T. Sferruzza, A. Tomassini, Bott–Chern Formality and Massey Products on Strong
K\"ahler with Torsion and K\"ahler Solvmanifolds, {\em J. Geom. Anal.} {\bf 34} (2024), https://doi.org/10.1007/s12220-024-01764-w
\bibitem{STc} T. Sferruzza, A. Tomassini, Hermitian geometrically Bott-Chern formal manifolds, {\tt arXiv:2507.10263 [math.DG]} (2025). %\bibitem{Ste21} J. Stelzig, On the structure of double complexes, %\emph{J. London Math. Soc. (2)} {\bf 104} (2021) 956–-988. %doi:10.1112/jlms.12453. https://doi.org/10.1112/jlms.12453 	 
%\bibitem{SW} J. Stelzig, S.O. Wilson, A $dd^c$-type condition beyond the K\"ahler realm, $\mathtt{arXiv:2208.01074v1}$ [math.DG].
\bibitem{Sul75} D. Sullivan, Differential Forms and the Topology of Manifolds,  in \emph{Manifolds}, 
Tokyo 1973, ed. A. Hattori, Univ. Tokyo Press,  1975.
\bibitem{Sul77} D. Sullivan, Infinitesimal computations in topology, \emph{Publ. Mat. Inst. Hautes Études Sci.} {\bf 47} (1977), p. 269--331, https://doi.org/10.1007/BF02684341
\bibitem{TarTom} N. Tardini, A. Tomassini, On geometric Bott-Chern formality and deformations, {\em Ann. Mat. Pura ed Appl.} {\bf 196} (2017), 349--362, https://doi.org/10.1007/s10231-016-0575-6
%\bibitem{TomTor} A. Tomassini, S. Torelli, On Dolbeault formality %and small deformations, \emph{ Internat.  J. Math.} {\bf 25} %(2014), 1450111, https://doi.org/10.1142/S0129167X14501110
\end{thebibliography}
\end{document}